\documentclass[twoside,12pt]{article}

\def\ifundefined#1{\expandafter\ifx\csname#1\endcsname\relax}
\makeatletter
\RequirePackage{theorem}

\newcommand{\reportenum}[3][ ]{\gdef\rep@rtenum{#2}
\gdef\rep@rteyear{#3}\gdef\wh@reappear{#1}}

{\let\@ldmaketitle=\maketitle
\renewcommand{\maketitle}{{\def\newpage{}
{\scriptsize\parbox[t]{0.3\textwidth}{\noindent Reporte Interno
\# \rep@rtenum\\
Departamento de Matem\'aticas\\CINVESTAV del IPN\\Mexico City, \rep@rteyear}
\hfill
\parbox[t]{0.5\textwidth}{\wh@reappear}}
\@ldmaketitle}}
}
{\AtBeginDocument{
   \makeatletter
   \pagestyle{myheadings}
   \markboth{\hfill\ifundefined{@authorshort} \@author
                   \else \@authorshort
                   \fi\hfill}
            {\hfill\ifundefined{@titleshort} \@title
                   \else \@titleshort
                   \fi\hfill}
   \makeatother
 }
}

     \theorembodyfont{\slshape}
        \newtheorem{thm}{Theorem}[section]
     \newtheorem{prop}[thm]{Proposition}
     
     \newtheorem{cor}[thm]{Corollary}

     \theorembodyfont{\upshape}
     \newtheorem{defn}[thm]{Definition}
     
     \newtheorem{example}[thm]{Example}

     \newtheorem{rem}[thm]{\mdseries\scshape Remark}

\def\@addpunct#1{\ifnum\spacefactor>\@m \else#1\fi}

\newenvironment{proof}[1][\proofname]{\par
  \normalfont
  \topsep6\p@\@plus6\p@ \trivlist
  \item[\hskip\labelsep\scshape
    #1\@addpunct{.}]\ignorespaces
}{%
  $\qed$\endtrivlist
}
\newcommand{\proofname}{Proof}

\providecommand{\dedicatory}[1]{}
\providecommand{\keywords}[1]{\begingroup \def \protect {\noexpand \protect
\noexpand }\xdef \@thefnmark { }\endgroup \@footnotetext{{\em Keywords and
phrases.\/} #1}}
\providecommand{\AMSMSC}[2]{\begingroup \def \protect {\noexpand \protect
\noexpand }\xdef \@thefnmark { }\endgroup \@footnotetext{{1991 \it
Mathematical Subject Classification.\/} Primary: #1; Secondary: #2.}}

\newcommand{\authorshort}[1]{\gdef\@authorshort{#1}}
\newcommand{\titleshort}[1]{\gdef\@titleshort{#1}}

   {\@definecounter{equation}}{\@newctr {equation}[section]}

\newcommand{\comment}[1]{}
\newcommand{\smnegsk}{\mskip -2mu plus 1mu}
\usepackage{amsfonts}
\newcommand{\algebra}[1]{\ensuremath{{\mathfrak #1}}}

\newcommand{\Cliff}[2]{\ensuremath{{\bf Cl}(#1,#2)}}

\newcommand{\Space}[2]{\ensuremath{ {{\mathbb #1}^{#2}} }}

\ifundefined{qed}
    \DeclareMathSymbol{\qed}{0}{AMSa}{"03}
\fi
\newcommand{\fourier}[1]{\ensuremath{{\mathcal F}^{#1}}}

\newcommand{\norm}[1]{\parallel\smnegsk #1 \smnegsk\parallel}
\newcommand{\modulus}[1]{\mid\smnegsk #1 \smnegsk\mid}
\newcommand{\scalar}[2]{\langle #1,#2\rangle}

\providecommand{\eqref}[1]{\textup{(\ref{#1})}}

\authorshort{Vladimir V. Kisil}
\titleshort{Connection between Different Function Theories}
\reportenum{173}{1995}

\begin{document}
\title{Connection between Different Function Theories in Clifford
Analysis\thanks{This work was partially supported
by CONACYT
Project 1821-E9211, Mexico.}}

\author{Vladimir V. Kisil\thanks{On leave from the Odessa State
University.}\\
Departamento de Matem\'aticas,
CINVESTAV del I.P.N.,\\
Apartado Postal 14-740,
07000, M\'exico, D.F. M\'EXICO  \\
\large fax: (525)-752-6412\
e-mail: \ttfamily vkisil@mvax1.red.cinvestav.mx}

\date{January 2, 1995}
\maketitle
\begin{abstract}
We describe an explicit connection between solutions to equations
$Df=0$ (the Generalized Cauchy-Riemann equation) and $(D+M)f=0$, where
operators $D$ and
$M$ commute. The described connection allows to construct a ``function
theory'' (the Cauchy theorem, the Cauchy integral, the Taylor and
Laurent series etc.) for solutions of the second equation from the
known function theory for solution of the first
(generalized Cauchy-Riemann) equation.

As well known, many physical equations related to the orthogonal group
of rotations or the Lorentz group (the Dirac equation, the Maxwell
equation etc.) can be naturally formulated in terms of the Clifford
algebra. For them our approach gives an explicit connection between
solutions  with zero and non-zero mass (or external fields) and
provides with a family of formulas for calculations.
\AMSMSC{30G35}{34L40, 81Q05}
\end{abstract}
\newpage
\section{Introduction}

Clifford analysis~\cite{BraDelSom82,DelSomSou92} (the function theory
of nullsolutions of the generalized Cauchy-Riemann operator $D$) has a
structure closer that of
the complex analysis of one variable than to standard complex analysis
of several variables. The success of Clifford analysis is mainly
explained because
the generalized Cauchy-Riemann operator $D$~\eqref{eq:dirac}
factorizes the Laplace operator
$\Delta=\sum_0^n\frac{\partial ^2}{\partial x_j^2}$. The development
of Clifford analysis is motivated not only by mathematical reasons but
also by its high applicability to the various physic problem (the
Maxwell equations, the Dirac equation and other, which is connected
with the orthogonal group of rotations or the Lorentz group).
Applications of Clifford analysis to the operator theory may be found
at the~\cite{McInPryde87,KisRam95a}.
After the development of the function theory for nullsolutions of the
generalized Cauchy-Riemann operator it is natural to
look~\cite{BraSomVAck94,Zhenyuan91}
for an analogous function theory for $\lambda$-solutions of the
generalized Cauchy-Riemann
operator (or nullsolutions of the operator $D-\lambda I$), namely,
such a function $f$ that $Df=\lambda f$ (or $Df=f \lambda $ due to
non-commutativity of the Clifford multiplication).

In the paper we describe (Section~\ref{se:connection}) an explicit
connection between solutions to
equations $Df=0$ and $(D+M)f=0$, where operators $D$ and $M$ commute.
This particularly includes the mentioned above case of equation $Df=f
\lambda $. The described connection allows (Section~\ref{se:theory})
to construct a ``function
theory'' (the Cauchy theorem, the Cauchy integral, the Taylor and
Laurent series etc.) for solutions of the second equation from the
known function theory for solution of the first (generalized
Cauchy-Riemann) equation. For the mentioned above physical
applications our approach gives an explicit connection between
solutions  with zero and non-zero mass (or external fields) and
provides with a family of formulas for calculations.

We
would like to present here basic ideas rather than to achieve (even if
it is possible) the final level of generality. Thus many presented
results may be considered in another (sometimes wider) context. For
example, this gives an alternative approach to the theory of operator
$D_\alpha $ in quaternionic analysis~\cite{KraSha95a}.

Author gladly acknowledges helpful discussions of Clifford analysis
with V.~V.~Kravchenko, J.~Ryan, M.~V.~Shapiro, F.~Sommen, and
N.~L.~Vasilevski.

\section{A Connection between $M$-Solutions to the Generalized
Cauchy-Riemann Operator for
Different $M$}\label{se:connection}

Let now $e_j$ be generators of the Clifford algebra \Cliff{0}{n} (we
use
books~\cite{BraDelSom82,DelSomSou92} as a standard reference).
This
means that the following {\em anti-commutation\/} relations hold:
\begin{equation}\label{eq:anti-comm}
\{e_i,e_j\}:=e_i e_j+e_j e_i=-2\delta_{ij}e_0,
\end{equation}
where $e_0=I$. We will referee to $e_j$ as to operators (of orthogonal
transformations) acting in a (finite-dimensional) Hilbert
space~\cite[Chap.~0, \S~A.1]{DelSomSou92}.

Function $f:\Space{R}{n}\rightarrow\Cliff{0}{n}$ is
called {\em monogenic\/} if it satisfies the {\em generalized
Cauchy-Riemann equation\/}
\begin{equation}\label{eq:dirac}
Df:=\frac{\partial f(y)}{\partial y_0} -\sum_{j=1}^n e_j\frac{\partial
f(y)}{\partial y_j}=0\ \mbox{ or }\ \frac{\partial f(y)}{\partial y_0}
=\sum_{j=1}^n e_j\frac{\partial
f(y)}{\partial y_j}.
\end{equation}
The success of Clifford analysis is mainly explained because
the generalized Cauchy-Riemann operator~\eqref{eq:dirac} factorizes
the Laplace operator
$\Delta=\sum_0^n\frac{\partial ^2}{\partial x_j^2}$.

H\"ormander's remark from paper~\cite{Anderson69} gives us by
the
fundamental solution to the generalized Cauchy-Riemann equation in the
form
\begin{eqnarray*}
K(y)&=&\fourier{\eta\rightarrow y} e^{-iy_0\sum_{j=1}^n\eta_je_j}\\
    &=&\int_\Space{R}{n}e^{i\sum_{j=1}^n y_j\eta_j}\,
e^{-iy_0\sum_{j=1}^n\eta_je_j}\, d\eta.
\end{eqnarray*}
(``Simply take the Fourier transform with respect to the spatial
variables, and solve the equation in
$y_0$''~\cite{Anderson69}). Otherwise, any solution $f(y)$
to~\eqref{eq:dirac} is given by a
convolution of some function $\widetilde{f}(y)$ on \Space{R}{n-1}  and
the fundamental
solution $K(y)$. On the contrary, a convolution $K(y)$ with any
function
is a solution to~\eqref{eq:dirac}. We have:
\begin{eqnarray}
[K*f](y)&=&\int_\Space{R}{n}K(y-t)\,f(t)\,dt\nonumber \\
        &=&\int_\Space{R}{n}\int_\Space{R}{n}e^{-i\sum_{j=1}^n
(y_j-t_j)\eta_j}\, e^{-iy_0\sum_{j=1}^n\eta_je_j}\,
d\eta f(t)\,dt\nonumber\\
        &=&\int_\Space{R}{n}\int_\Space{R}{n}e^{i\sum_{j=1}^n
t_j\eta_j}\,
e^{-i\sum_{j=1}^n\eta_j(y_0 e_j-y_j e_0)}\, d\eta f(t)\,dt\nonumber \\
        &=&\int_\Space{R}{n} e^{-i\sum_{j=1}^n\eta_j(y_0 e_j-y_j
e_0)}\int_\Space{R}{n}e^{i\sum_{j=1}^n t_j\eta_j}\, f(t)\,dt\, d\eta
\nonumber \\
        &=&\int_\Space{R}{n} e^{-i\sum_{j=1}^n\eta_j(y_0 e_j-y_j
e_0)}\,
\widehat{f}(-\eta)\, d\eta. \label{eq:malonek}
\end{eqnarray}
Equation~\eqref{eq:malonek} defines the Weyl functional
calculus~\cite{Anderson69} for the function  $f(-y)$ and the $n$-tuple
of operators
\begin{equation}\label{eq:monom}
\{\vec{y}_j=y_0 e_j-y_j e_0\},\ 1\leq j \leq n.
\end{equation}
Thus any solution to the generalized Cauchy-Riemann
equation~\eqref{eq:dirac} can be written as a function of $n$
monomials~\eqref{eq:monom}:
\begin{equation}\label{eq:laville}
\breve{f}(y_0,y_1,\ldots,y_n)= f(y_0 e_1-y_1 e_0, y_0 e_2-y_2
e_0,\ldots,y_0 e_n-y_n e_0).
\end{equation}
Another significant remark: if we fix the value $y_0=0$
in~\eqref{eq:laville} we easily obtain:
\begin{displaymath}
\breve{f}(0,y_1,\ldots,y_n)= f(-y_1 e_0, -y_2 e_0,\ldots,-y_n e_0)=
f(-y_1 , -y_2 ,\ldots,-y_n )e_0.
\end{displaymath}
Thus we may consider the function $\breve{f}(y_0,y_1,\ldots,y_n)$ in
$n+1$ variables as {\em analytic (or Cauchy-Kovalevska) expansion\/}
for the function $f(y_1,\ldots,y_n)$ in $n$ variables (compare
with~\cite{Laville91}).

Using the power series decomposition for the exponent one can see that
formula~\eqref{eq:malonek} defines the permutational (symmetric)
product
of monomials~\eqref{eq:monom}. The significant role of such monomials
and functions generated by them were described for  quaternionic
analysis
in~\cite{Sudbery79}, for Clifford analysis
in~\cite{Laville87,Malonek93},
for Fueter-Hurwitz analysis in~\cite{KrolRam92}.
 But during our consideration we
used only the commutation relation $[e_0,e_j]=0$ and never used the
anti-commutation relations~\eqref{eq:anti-comm}. Thus
formula~\eqref{eq:malonek} {\em is true and may be useful without
Clifford analysis\/} as well.
\begin{prop}
Any solution to equation~\eqref{eq:dirac}, where $e_j$ are arbitrary
self-adjoint operators, is given as arbitrary function of $n$
monomials~\eqref{eq:monom} by the formula~\eqref{eq:malonek}.
\end{prop}

Due to physical application we will consider equation
\begin{equation}\label{eq:mass}
\frac{\partial f}{\partial y_0}=(\sum_{j=1}^n e_j\frac{\partial
}{\partial y_j}+M)f,
\end{equation}
where $e_j$ are arbitrary self-adjoint operators and $M$ is a bounded
operator commuting with all $e_j\frac{\partial}{\partial y_j}$.
\begin{example}
If $e_j$ are generators~\eqref{eq:anti-comm} of the Clifford algebra
and $M=M_\lambda $ is an operator of
multiplication from the {\em right\/}-hand side by the Clifford number
$\lambda $,
differential operator
\begin{equation}\label{eq:D_M}
(\sum_{j=1}^n e_j\frac{\partial}{\partial y_j}+M)f,
\end{equation}
 factorizes the Helmholtz
operator
$\Delta +M_{\lambda^2}$. Equation~\eqref{eq:mass} is known in quantum
mechanics as the
{\em Dirac
equation for a particle with a non-zero rest
mass\/}~\cite[\S 20]{BerLif82}, \cite[\S 6.3]{BogShir80} and
\cite{Kravchenko95a}. We will specialize our results for the case
$M=M_\lambda$, especially for the simplest (but still important!) case
$\lambda\in\Space{R}{}$.
\end{example}

\begin{example}
It is well known~\cite[Chap.~VI, \S~1.2]{KirGvi88}, that an operator
commutes with differential operators (particularly with the
generalized Cauchy-Riemann
operator) if and only if it is an operator of convolution with a
generalized function. The operator $M$ can be used for an introducing
into equation~\eqref{eq:mass} a non-zero mass of the particle or
external fields.

Notably, the operator $M_\lambda$ from the
previous Example is the convolution from the right-hand side
(non-commutativity of Clifford multiplication!) with the Dirac
function $\lambda \delta (x)$.
\end{example}

Thus the family of possible operators $M$ is rather width and includes
examples with interesting applications.
Simple modification of the previous calculations~\eqref{eq:malonek}
gives us the following result
\begin{prop}
Any solution to equation~\eqref{eq:mass}, where $e_j$ are arbitrary
self-adjoint operators and $M$ commutes with all
$e_j\frac{\partial}{\partial y_j}$, is given by the
formula
\begin{equation}\label{eq:m}
e^{y_0 M}\int_\Space{R}{n} e^{-i\sum_{j=1}^n\eta_j(y_0 e_j-y_j e_0)}
\widehat{f}(-\eta)\, d\eta,
\end{equation}
where $f$ is an arbitrary function on \Space{R}{n-1}.
\end{prop}

Here
\begin{equation}\label{eq:exponent}
e^{y_0 M}=\sum_{j=0}^{\infty} \frac{(y_0 M)^j}{j!}
\end{equation}
is well defined for all bounded $M$. Comparing~\eqref{eq:malonek}
and~\eqref{eq:m} we obtain
\begin{thm}
The function $f(y)$ is a solution to the equation
\begin{displaymath}
\frac{\partial f}{\partial y_0}=(\sum_{j=1}^n e_j\frac{\partial
}{\partial y_j}+M_1)f
\end{displaymath}
 if and only if the function
\begin{displaymath}
g(y)=e^{y_0 M_2} e^{-y_0 M_1} f(y)
\end{displaymath}
is a solution to the equation
\begin{displaymath}
\frac{\partial g}{\partial y_0}=(\sum_{j=1}^n e_j\frac{\partial
}{\partial y_j}+M_2)g,
\end{displaymath}
where $M_1$ and $M_2$ are bounded operators commuting with $e_j$.
\end{thm}
\begin{cor}\label{co:mass}
The function $f(y)$ is a solution to the equation~\eqref{eq:mass}
 if and only if the function
$e^{y_0 M}f(y)$
is a solution to the generalized Cauchy-Riemann
equation~\eqref{eq:dirac}.

In the case $M=M_\lambda$ we have $e^{y_0 M_\lambda}f(y)=f(y)e^{y_0
\lambda}$ and if $\lambda\in \Space{R}{}$ then $e^{y_0
M_\lambda}f(y)=f(y)e^{y_0 \lambda}=e^{y_0 \lambda}f(y)$.
\end{cor}

\section{Function Theory for $M$-Solutions of the Generalized
Cauchy-Riemann Operator}\label{se:theory}

In this section we construct a function theory (in the classic sense)
for $M$-solutions of
the generalized Cauchy-Riemann operator basing on Clifford analysis
and
Corollary~\ref{co:mass}. Proofs are very short and almost evident, but
we present them for the sake of completeness.

The set of solutions to~\eqref{eq:dirac} and~\eqref{eq:mass} in a nice
domain $\Omega$ will be denoted by
$\algebra{M}(\Omega)=\algebra{M}_0(\Omega)$ and
$\algebra{M}_M(\Omega)$ correspondingly. In the case $M=M_\lambda$ we
use the notation
$\algebra{M}_\lambda(\Omega)=\algebra{M}_{M_\lambda}(\Omega)$ also. We
suppose that all functions from $\algebra{M}_\lambda(\Omega)$ are
continuous in the closure of $\Omega$. Let
\begin{equation}\label{eq:cauchy-ker}
E(y-x)=
\frac{\Gamma(\frac{n+1}{2})}{2\pi^{(m+1)/2}}\,
\frac{\overline{y-x}}{\modulus{y-x}^{n+1}}
\end{equation}
be the Cauchy kernel~\cite[p.~146]{DelSomSou92}  and
\begin{equation}
d\sigma=\sum_{j=0}^n (-1)^j e_j dx_0 \wedge \ldots \wedge[dx_j] \wedge
\ldots \wedge dx_m.
\end{equation}
be the differential form of the ``oriented surface
element''~\cite[p.~144]{DelSomSou92}. Then for any
$f(x)\in\algebra{M}(\Omega)$ we have the Cauchy integral
formula~\cite[p.~147]{DelSomSou92}
\begin{equation}
\int_{\partial \Omega} E(y-x)\,d\sigma_y\,
f(y)=\left\{\begin{array}{cl}
f(x),& x\in\Omega\\
0,& x\not\in\bar{\Omega}
\end{array}.\right.
\end{equation}

\begin{thm}[Cauchy's Theorem]
Let $f(y)\in \algebra{M}_M(\Omega)$. Then
\begin{displaymath}
\int_{\partial \Omega} d\sigma_y\,e^{-y_0 M}f(y)=0.
\end{displaymath}
Particularly, for $f(y)\in \algebra{M}_\lambda(\Omega)$ we have
\begin{displaymath}
\int_{\partial \Omega} d\sigma_y\,f(y)e^{-y_0 \lambda}=0,
\end{displaymath}
and
\begin{displaymath}
\int_{\partial \Omega} d\sigma_ye^{-y_0 \lambda}\,f(y)=0,
\end{displaymath}
if $\lambda\in\Space{R}{} $.
\end{thm}

\begin{proof}
It easily follows because $e^{-y_0 M}f(y)\in
\algebra{M}(\Omega)$ and the corresponding result for the generalized
Cauchy-Riemann
equation~\cite[Chap.~II, \S~0.2.1]{DelSomSou92}.
\end{proof}

\begin{thm}[Cauchy's Integral Formula]
Let $f(y)\in \algebra{M}_M(\Omega)$. Then
\begin{equation}\label{eq:m-cauchy}
e^{x_0 M}\int_{\partial \Omega} E(y-x)\,d\sigma_y\,
e^{-y_0 M}f(y)=\left\{\begin{array}{cl}
f(x),& x\in\Omega\\
0,& x\not\in\bar{\Omega}
\end{array}.\right.
\end{equation}
Particularly, for $f(y)\in \algebra{M}_\lambda(\Omega)$ we have
\begin{displaymath}
\int_{\partial \Omega} E(y-x)\,d\sigma_y\,
f(y)e^{(x_0-y_0) \lambda}=\left\{\begin{array}{cl}
f(x),& x\in\Omega\\
0, &x\not\in\bar{\Omega}
\end{array}.\right.
\end{displaymath}
and
\begin{displaymath}
\int_{\partial \Omega} E(y-x)e^{(x_0-y_0) \lambda}\,d\sigma_y\,
f(y)=\left\{\begin{array}{cl}
f(x),& x\in\Omega\\
0,& x\not\in\bar{\Omega}
\end{array}.\right.
\end{displaymath}
if $\lambda\in\Space{R}{} $.
\end{thm}
\begin{proof}
It easily follows from the fact that $e^{-y_0 M}f(y)\in
\algebra{M}(\Omega)$ and the corresponding result for the generalized
Cauchy-Riemann
equation~\cite[Chap.~II, \S~0.2.2]{DelSomSou92}.
\end{proof}
It is hard to expect that formula~\eqref{eq:m-cauchy} may be rewritten
as
\begin{displaymath}
\int_{\partial \Omega} E'(y-x)\,d\sigma_y\,
f(y)=\left\{\begin{array}{cl}
f(x),& x\in\Omega\\
0,& x\not\in\bar{\Omega}
\end{array}\right.
\end{displaymath}
with a simple function $E'(y-x)$.

Because an application of the bounded operator $e^{y_0 M}$ does not
destroy uniformed convergency of functions we obtain
(cf.~\cite[Chap.~II, \S~0.2.2, Theorem~2]{DelSomSou92})
\begin{thm}[Weierstrass' Theorem]
Let $\{f_k\}_{k\in\Space{N}{}}$ be a sequence in
$\algebra{M}_M(\Omega)$, which converges uniformly to $f$ on each
compact subset $K\in \Omega$. Then
\begin{enumerate}
\item $f\in \algebra{M}_M(\Omega)$.
\item For each multi-index
$\beta=(\beta_0,\ldots,\beta_m)\in\Space{N}{n+1}$, the sequence
$\{\partial ^\beta f_k\}_{k\in\Space{N}{}} $ converges uniformly on
each compact subset $K\in \Omega$ to $\partial ^\beta f$.
\end{enumerate}
\end{thm}

\begin{thm}[Mean Value Theorem]
Let $f\in\algebra{M}_M(\Omega)$. Then for all $x\in \Omega$ and $R>0$
such that the ball $\Space{B}{}(x,R)\in\Omega $,
\begin{displaymath}
f(x)= e^{x_0 M}
\frac{(n+1)\Gamma(\frac{n+1}{2})}{2R^{n+1}\pi^{(m+1)/2}}
\int_{\Space{B}{}(x,R)} e^{-y_0 M} f(y)\,dy.
\end{displaymath}
\end{thm}
\begin{proof}
We again referee to corresponding theorem~~\cite[Chap.~II, \S~0.2.2,
Theorem~3]{DelSomSou92} for $\algebra{M}(\Omega)$ and
Corollary~\ref{co:mass}.
\end{proof}

We skip Morera's, Painlev\'e's theorems~\cite[Chap.~II,
\S~0.2.3]{DelSomSou92} and Plemelj-Sokhotskij formulas formulated for
$\algebra{M}_M(\Omega)$.

It is possible also to introduce the notion of
the {\em differentiability\/}~\cite{Malonek93}  for
solutions to~\eqref{eq:dirac}, namely, an increment of any solution
to~\eqref{eq:dirac} may be approximated up to infinitesimals of the
second
order by a linear function of monomials $\vec{y}_j$
from~\eqref{eq:monom}. We now give $\lambda$-version\footnote{At the
rest of the paper we will consider only the case $M=M_\lambda$. The
general case needs only obvious modifications.} of this notion.

\begin{defn}
The function $f(y)$ is called {\em $\lambda$-differentiable\/} from
the left at $y$ if there exist a $\lambda$-linear form
\begin{equation}
l(\vec{y})=\sum_{j=1}^n \vec{y}_j A_j e^{y_0 \lambda},\
A_j\in\Cliff{n}{0}
\end{equation}
such that
\begin{equation}
\lim_{\Delta\vec{y}\rightarrow 0 }
\frac{\modulus{f(\vec{y}-\Delta\vec{y})-
f(\vec{y})- l(\vec{y})}}{\norm{\Delta\vec{y}}}=0.
\end{equation}
The function is $\lambda$-differentiable in $\Omega$ if the function
is $\lambda$-differentiable at all points of $\Omega$. The
$\lambda$-linear mapping $l(\vec{y})$ is called the
{\em $\lambda$-derivative\/} of the function $f$ at the point $y$.
\end{defn}

We skip proofs of the following main results connected with the notion
of $\lambda$-differentiability (compare with~\cite{Malonek93}).

\begin{thm}
If $f(\vec{y})$ is $\lambda$-differentiable then the corresponding
$\lambda$-derivative is determined in a unique way.
\end{thm}

\begin{thm}
Let $f(\vec{y})$ be continuously real differentiable at a point
$\vec{y}$. The function $f$ is $\lambda$-differentiable at $\vec{y}$
if and only if $f\in\algebra{M}_\lambda$.
\end{thm}

Let us remind that symmetric
products~{\cite{Malonek93}} is defined by the formula
\begin{equation}\label{eq:permut}
a_1\times a_2\times \cdots \times a_k=\frac{1}{k!}\sum a_{j_1} a_{j_2}
\cdots  a_{j_k},
\end{equation}
where the sum is taken over all of permutations of $({j_1} ,{j_2},
\ldots,
{j_k})$. For a multi-index $\beta$ the notation $\vec{y}^\beta$
denotes the corresponding symmetric product. If we consider the linear
combination of such products then Clifford valued coefficients are
written on the right-hand
side.

\begin{thm}
Every $f(\vec{y})\in\algebra{M}_\lambda(\Omega)$ is infinitely
$\lambda$-differentiable and for any point $a\in\Omega$ can be
presented by the Taylor series
\begin{displaymath}
f(\vec{y})=\sum_{\beta=0}^\infty (\vec{y}-\vec{a})^\beta c_\beta
e^{-y_0 \lambda}
\end{displaymath}
in some neighborhood of $a$.
\end{thm}

An important property of Clifford analysis is the existence of the
reproducing Bergman kernel~\cite[\S~24]{BraDelSom82}. We will give the
$\lambda$-version of this
result (compare with~\cite{BraSomVAck94,ShaVas93a}).
\begin{thm}
$\algebra{M}_\lambda(\Omega)$ has the reproducing formula
\begin{equation}\label{eq:m-bergman}
\int_\Omega B(x,y) \, f(y)e^{(x_0-y_0) \lambda}\, dy=f(x),
\end{equation}
where $B(x,y)$ is the Bergman kernel~\cite[\S~24.1]{BraDelSom82} from
Clifford analysis.

Moreover, if $\lambda\in\Space{R}{}$ then formula~\eqref{eq:m-bergman}
takes the form of usual reproducing formula
\begin{equation}
\int_\Omega B'_\lambda(x,y) \, f(y)\, dy=f(x),
\end{equation}
where $B'_\lambda(x,y)=B(x,y)e^{(x_0-y_0) \lambda}$. For the unit ball
$\Space{B}{}(0,1)$ centered at the origin the explicit formula for
$B'_\lambda(x,y)$ is (compare with~\cite{ShaVas93a}):
\begin{eqnarray*}
B'_\lambda(x,y)&=&
\frac{\Gamma(\frac{n+1}{2}) (n+1) e^{(x_0-
y_0)\lambda}}{2\pi^{(n+1)/2}}\left( \frac{n+1}{(1-
2\scalar{y}{x}+\modulus{y}^2\cdot\modulus{x}^2)^{(n+1)/2}}\right.\\
&&-\frac{2\bar{x}y}
{(1-2\scalar{y}{x}+\modulus{y}^2\cdot\modulus{x}^2)^{(n+1)/2}}\\
&&\left. +\frac{(n+1)\overline{(y-x\modulus{y}^2)}(x-
y\modulus{x}^2)}{(1-
2\scalar{y}{x}+\modulus{y}^2\cdot\modulus{x}^2)^{(n+3)/2}}\right).
\end{eqnarray*}
\end{thm}

\begin{rem}
Quaternionic analysis is not (formally speaking) a corollary of
Clifford analysis. Thus quaternionic analysis needs an application the
ideas of this paper rather than concrete results.
\end{rem}
\newcommand{\noopsort}[1]{} \newcommand{\printfirst}[2]{#1}
  \newcommand{\singleletter}[1]{#1} \newcommand{\switchargs}[2]{#2#1}
  \newcommand{\irm}{\mbox{\rm I}} \newcommand{\iirm}{\mbox{\rm II}}
  \newcommand{\vrm}{\mbox{\rm V}}

\end{document}